\documentclass[10pt]{article}
\usepackage{graphicx}
\usepackage{amsthm,amsmath,amsfonts,amssymb} % added by PB
\usepackage{color}
\catcode`\@=11 \@addtoreset{equation}{section}

\catcode`\@=12

\newtheorem{Theorem}{Theorem}[section]
\newtheorem{Definition}{Definition}[section]
\newtheorem{Proposition}{Proposition}[section]
\newtheorem{Lemma}{Lemma}[section]
\newtheorem{Corollary}{Corollary}[section]
\newtheorem{Remark}{Remark}[section]

\newcommand{\bTheorem}[1]{
\begin{Theorem} \label{T#1} }
\newcommand{\eT}{\end{Theorem}}

\newcommand{\bDefinition}[1]{
\begin{Definition} \label{D#1} }
\newcommand{\eD}{\end{Definition}}

\newcommand{\bProposition}[1]{
\begin{Proposition} \label{P#1}}
\newcommand{\eP}{\end{Proposition}}

\newcommand{\bRemark}[1]{
\begin{Remark} \label{R#1}}
\newcommand{\eR}{\end{Remark}}

\newcommand{\bLemma}[1]{
\begin{Lemma} \label{L#1} }
\newcommand{\eL}{\end{Lemma}}

\newcommand{\bCorollary}[1]{
\begin{Corollary} \label{C#1} }
\newcommand{\eC}{\end{Corollary}}

\newcommand{\bFormula}[1]{
\begin{equation} \label{#1}}
\newcommand{\eF}{\end{equation}}

\newcommand{\Ov}[1]{\overline{#1}}

\newcommand{\DC}{C^\infty_c}

\newcommand{\vr}{\varrho}

\newcommand{\vL}{\vc{L}}
\newcommand{\vu}{\vc{u}}
\newcommand{\vc}[1]{{\bf #1}}
\newcommand{\Div}{{\rm div}_x}
\newcommand{\Grad}{\nabla_x}

\newcommand{\tn}[1]{\mbox {\F #1}}
\newcommand{\dx}{{\rm d} {x}}
\newcommand{\dt}{{\rm d} t }

\newcommand{\vJ}{\vc{J}}

\newcommand{\dxdt}{\dx \ \dt}
\newcommand{\intO}[1]{\int_{\Omega} #1 \ \dx}

\newcommand{\ep}{\varepsilon}

\font\F=msbm10 scaled 1000
\newcommand{\R}{\mathbb{R}}
\newcommand{\Z}{\mathbb{Z}}

\newcommand{\Del}{\Delta_x}

\definecolor{Cgrey}{rgb}{0.85,0.85,0.85}
\definecolor{Cblue}{rgb}{0.50,0.85,0.85}
\definecolor{Cred}{rgb}{1,0,0}
\definecolor{fancy}{rgb}{0.10,0.85,0.10}

\newcommand\Cbox[2]{%
    \newbox\contentbox%
    \newbox\bkgdbox%
    \setbox\contentbox\hbox to \hsize{%
        \vtop{
            \kern\columnsep
            \hbox to \hsize{%
                \kern\columnsep%
                \advance\hsize by -2\columnsep%
                \setlength{\textwidth}{\hsize}%
                \vbox{
                    \parskip=\baselineskip
                    \parindent=0bp
                    #2}%
                \kern\columnsep%
            }%
            \kern\columnsep%
        }%
    }%
    \setbox\bkgdbox\vbox{
        \color{#1}
        \hrule width  \wd\contentbox %
               height \ht\contentbox %
               depth  \dp\contentbox
        \color{black}
    }%
    \wd\bkgdbox=0bp%
    \vbox{\hbox to \hsize{\box\bkgdbox\box\contentbox}}%
    \vskip\baselineskip%
}

\date{}

\begin{document}

%%%%%%%%%%%%%%%%%%%%%%%%%%%%%%%%

\title{Well/ill posedness for the Euler-Korteweg-Poisson system and related problems}

\author{Donatella Donatelli \and
Eduard Feireisl \thanks{The research of E.F. leading to these results has received funding from the European Research Council under the European Union's Seventh Framework
Programme (FP7/2007-2013)/ ERC Grant Agreement 320078. The Institute of Mathematics of the Academy of Sciences of the Czech
Republic is supported by RVO:67985840.} 
\and Pierangelo Marcati \thanks{The research of D.D. and P.M leading to these results was partially supported bt the Grant PRIN N. 2012LWXHJ}}

%\thanks{Eduard Feireisl acknowledges the support of the GA\v CR (Czech Science Foundation) project P201-13-00522S in the framework of RVO: %67985840.} }

%\thanks{Eduard Feireisl acknowledges the support of the project LL1202 in the
%programme ERC-CZ funded
%by the Ministry of Education, Youth and Sports of the Czech Republic.}

%\thanks{The research of E.F. leading to these results has received funding from the European Research Council under the European Union's %Seventh Framework
%Programme (FP7/2007-2013)/ ERC Grant Agreement 320078}

\maketitle

\centerline{Department of Information Engineering, Computer Science and Mathematics}
\centerline{University of L'Aquila, 67100 L'Aquila, Italy}

\bigskip

\centerline{and}
\centerline{Institute of Mathematics of the Academy of Sciences of the Czech Republic}
\centerline{\v Zitn\' a 25, 115 67 Praha 1, Czech Republic}

%\centerline{Charles University in Prague, Faculty of Mathematics and Physics, Mathematical Institute}

%\centerline{Sokolovsk\' a 83, 186 75 Praha 8, Czech Republic}

\bigskip

\centerline{and}

\centerline{Department of Information Engineering, Computer Science and Mathematics}
\centerline{University of L'Aquila}
\centerline{and GSSI - Gran Sasso Science Institute, 67100 L'Aquila, Italy}

\bigskip

\begin{abstract}

We consider a general Euler-Korteweg-Poisson system in $R^3$, supplemented with the space periodic boundary conditions, where the quantum hydrodynamics equations and the classical fluid dynamics equations with capillarity are recovered as particular examples.  We show that the system admits infinitely many global-in-time weak solutions for any sufficiently smooth
initial data including the case of a vanishing initial density - the vacuum zones. Moreover, there is a vast family of initial data, for which the Cauchy problem possesses infinitely many dissipative weak solutions, i.e. the
weak solutions satisfying the energy inequality. Finally, we establish the weak-strong uniqueness property in a class of solutions without vacuum. In this paper we show that, even in presence of a dispersive tensor, we have the same phenomena found by De Lellis and Sz{\'e}kelyhidi.

\end{abstract}

\medskip

%\subjclass[2010]{Primary: 35Q30, 76N10; Secondary: 46E35.}

\medskip

%\tableofcontents

\section{Introduction}
\label{i}

A general \textsc{Euler-Korteweg-Poisson system} describing the time evolution of the density $\vr = \vr(t,x)$ and the momentum $\vc{J} = \vc{J}(t,x)$ of
an inviscid fluid can be written in the form:

%\Cbox{Cgrey}{
\bFormula{i1}
\partial_t \vr + \Div \vc{J} = 0,
\eF
\bFormula{i2}
\partial_t \vc{J} + \Div \left( \frac{\vc{J} \times \vc{J} }{\vr} \right) + \Grad p(\vr) =  - \alpha \vc{J} +  \vr \Grad \left(
K(\vr) \Del \vr + \frac{1}{2} K'(\vr) |\Grad \vr |^2 \right) + \vr \Grad V,
\eF
\bFormula{i3}
\Del V = \vr - \Ov{\vr},
\eF
%}

\noindent
where $K : (0, \infty) \to (0, \infty)$ is a smooth function, see Audiard \cite{Audi1}, \cite{Audi2}, Benzoni-Gavage et al. \cite{Benz1}, \cite{Benz2}. In particular,
taking $K = \Ov{K} > 0$ a positive constant, we recover the standard equations of an \emph{inviscid capillary fluid} (see Bresch et al.
\cite{BDD}, Kotchote \cite{Kot1}, \cite{Kot2}), while the choice $K(\vr) = \frac{\hbar}{4 \vr}$ gives rise to the so-called \emph{quantum fluid system} (see for instance Antonelli and Marcati \cite{AntMar1}, \cite{AntMar2}, J\" ungel \cite[Chapter 14]{Jun} and the references therein). In the latter case, the equations (\ref{i1} - \ref{i3}), by using the Madelung transformations,
may be formally seen as a description of the evolution of the momenta
\bFormula{i3+}
\vr = |\psi|^2, \ \vc{J} = \hbar \Im [\Ov{\psi} \Grad \psi ],
\eF
where the wave function $\psi$, in the case $\alpha=0$ and $\Ov{\vr} = 0$, is a solution of the following Schr\" odinger-Poisson
system:
\bFormula{i3++}
\imath \hbar \partial_t \psi = - \frac{\hbar^2}{2} \Del \psi - V \psi + f(|\psi|^2)\psi ,
\ \Delta V = |\psi|^2,
\eF
provided $p'(\vr) = \vr f'(\vr)$.

For the sake of simplicity, we consider the system (\ref{i1} - \ref{i3}) supplemented with the spatially periodic boundary conditions, namely on
the ``flat'' torus

%\Cbox{Cgrey}{
\[
\Omega = \tn{T}^3 \equiv   {\R^3/ \Z^3},
\]
and with the initial state
\bFormula{i4}
\vr(0, \cdot) = \vr_0, \ \vc{J}(0, \cdot) = \vc{J}_0.
\eF
%}

In view of the applications to the quantum fluid models, we consider a general non-negative distribution of the density $\varrho$ including the
vacuum zones where $\vr = 0$. We note that the Korteweg tensor can be written in the form
\[
\vr \Grad \left(
K(\vr) \Del \vr + \frac{1}{2} K'(\vr) |\Grad \vr |^2 \right)
\]
\[
=
\Div \left[ \vr \Div \Big( K(\vr) \Grad \vr \Big) \tn{I} \right] +
\frac{1}{2} \Div \left[ \Big( K(\vr) - \vr K'(\vr) \Big) |\Grad \vr|^2 \tn{I} \right]
 - \Div \Big[
K(\vr) \Grad \vr \otimes \Grad \vr \Big].
\]
Thus, introducing
\bFormula{chi}
\chi(\vr) = \vr K(\vr),
\eF
we deduce that
\bFormula{i5}
\vr \Grad \left(
K(\vr) \Del \vr + \frac{1}{2} K'(\vr) |\Grad \vr |^2 \right)
\eF
\[
= \Grad \Big( \chi (\vr) \Del \vr  \Big) + \frac{1}{2} \Grad \Big( \chi'(\vr) |\Grad \vr|^2 \Big) - 4 \Div \Big( \chi(\vr) \Grad \sqrt{\vr} \otimes \Grad \sqrt{\vr} \Big) \equiv \Div \mathcal{K} (\vr, \Grad \vr),
\]
\bFormula{i6}
\mathcal{K}(\vr, \Grad \vr) = \Big[ \chi (\vr) \Del \vr + \frac{1}{2} \chi'(\vr) |\Grad \vr|^2 \Big] \tn{I} - 4 \chi(\vr) \Grad \sqrt{\vr} \otimes \Grad \sqrt{\vr},
\eF
where the choice $\chi \equiv \hbar /4 $ determines the quantum fluids while $\chi(\vr) = \vr$ corresponds to the capillary fluids with constant capillarity. Accordingly  the choice of $\chi(\vr) $ determines the role of the quadratic nonlinearities, in the case of the quantum fluids the term sensitive to the appearance of the vacuum, beyond the convective term, is then $\Grad \sqrt{\vr} \otimes \Grad \sqrt{\vr}$.

The parameter $\alpha \geq 0$ in (\ref{i2}) represent a damping effect relevant in certain applications, in particular collision effects for quantum models for semiconductor devices. In what follows, we assume, for the sake of simplicity, that $\alpha = 1$. Strangely enough, the presence of damping makes the problem more difficult in view of the methods used in the present paper and due to the dispersive  nature of the equations. We remark that the theory we develop below applies to the case $\alpha = 0$ as well, with only obvious modifications in the proofs.

\subsection{Energy}

The Euler-Korteweg-Poisson system (\ref{i1}-\ref{i3}) admits a natural energy density, namely

%\Cbox{Cgrey}{

\bFormula{i7}
E (\vr, \Grad \vr, \vc{J} ) = \frac{1}{2} \frac{ |\vc{J}|^2 }{\vr} + P(\vr) + \frac{ K(\vr) }{2}| \Grad \vr |^2  + \frac{1}{2} |\Grad V|^2
\eF
\[
= \frac{1}{2} \frac{ |\vc{J}|^2 }{\vr} + P(\vr) + 2 \chi(\vr)| \Grad \sqrt{\vr} |^2  + \frac{1}{2} |\Grad V|^2
\]

%}

\noindent
where $\chi$ was introduced in (\ref{chi}) and
\[
P(\vr) = \vr \int_{1}^\vr \frac{p(z)}{z^2} \ {\rm d}z.
\]
Indeed, taking the scalar product of the momentum equation (\ref{i2}) with $\vc{J}/ \vr$ and using (\ref{i1}), (\ref{i3}), we obtain the energy balance
\bFormula{i8}
\frac{{\rm d}}{{\rm d}t} \intO{ E (\vr, \Grad \vr, \vc{J} ) (t, \cdot) } + \intO{ \frac{|\vc{J} |^2 }{\vr}(t, \cdot) } = 0.
\eF

In this paper, we focus on \emph{bounded energy} (weak) solutions for which $E(\vr, \Grad \vr, \vJ)$ is bounded on the whole physical space $\Omega$ and for any time $t \in [0,T]$. In particular, the momentum $\vJ$ must vanish on the vacuum set where $\vr  = 0$.

\subsection{Velocity Fields}

As already pointed out several times, our goal is to consider the solutions that may contain vacuum zones. In the context of quantum hydrodynamics, the classical WKB formalism does not  allow the definition of the velocity in the nodal regions, while  the current measure  $\vc{J}dx$ obtained via the Madelung transform can be differentiated in the sense of measure in $\vr dx$ but the velocity field defined in this way is  ${\text L}^1(\vr dx)$ only. In the context of classical fluid mechanics, where vacuum is not permitted in the natural framework of applications of the model, it is customary to replace the momentum $\vc{J}$ by $\vr \vu$, where
$\vu$ is the macroscopic \emph{velocity} of the fluid. We emphasize that the velocity $\vu$ has a physical interpretation only on the sets where $\vr > 0$ and, in particular, it has no particular meaning on the vacuum. For these reasons, we avoid using the concept of velocity in the formulation of our problem and we are going to develop a self consistent theory in the $(\vr,\vc{J})$ variables.  The vacuum problem has been extensively discussed in \cite{AntMar1}, \cite{AntMar2}.

\subsection{Weak solutions}

Since the solutions of the problem (\ref{i1}-\ref{i3}), (\ref{i4}) may not be regular on the vacuum, quantum vortices may appear and moreover the hydrodynamic variables $(\vr,\vc{J})$ may not have better regularity than the energy space, it seems natural to introduce the concept of
\emph{weak solution}.

%\Cbox{Cgrey}{

\bDefinition{i1}

We say that
\[
\vr \in C_{\rm weak}([0,T]; L^2(\Omega)) \cap L^\infty((0,T) \times \Omega),\
\vJ \in C_{\rm weak}([0,T]; L^2(\Omega;R^3)) \cap L^\infty((0,T) \times \Omega ; R^3)
\]
is a \emph{bounded energy weak solution} to the Euler-Korteweg-Poisson system (\ref{i1}-\ref{i3}), (\ref{i4}) if
\bFormula{i9-}
\vr(t, \cdot) > 0 \ \mbox{a.a. in} \ \Omega \ \mbox{for any}\ t \in (0,T), \ \Grad \vr \in L^\infty((0,T) \times \Omega),
\eF
\[
E (\vr, \Grad \vr, \vJ)(t, \cdot) \leq \Ov{E} \ \mbox{for a.a.} \ t \in (0,T),
\]
\bFormula{i9+}
\vr(0, \cdot) = \vr_0, \ \vJ(0, \cdot) = \vJ_0,
\eF
and the following integral identities
\bFormula{i9}
- \intO{ \vr \varphi } \Big|^{t = \tau_2}_{t = \tau_1} +
\int_{\tau_1}^{\tau_2} \intO{ \left( \vr \partial_t \varphi + \vJ \cdot \Grad \varphi \right) } \ \dt = 0,
\eF
\bFormula{i10}
- \intO{ \vJ \cdot \varphi } \Big|^{t = \tau_2}_{t = \tau_1} +
\int_{\tau_1}^{\tau_2} \intO{ \left( \vJ \cdot \partial_t \varphi + \frac{ \vJ \otimes \vJ }{\vr} : \Grad \varphi + p(\vr) \Div \varphi \right) } \ \dt
\eF
\[
= \int_{\tau_1}^{\tau_2} \intO{ \left(  - \Grad \vr \cdot \Grad ( \chi (\vr)  \Div \varphi ) + \frac{1}{2} \chi'(\vr) |\Grad \vr|^2 \Div \varphi
 - 4 \chi(\vr) \Grad \sqrt{\vr} \otimes \Grad \sqrt{\vr} : \Grad \varphi \right) }
\]
\[
+ \int_{\tau_1}^{\tau_2} \intO{ \Big( \vc{J} \cdot \varphi - \vr \Grad \vc{V} \cdot \varphi \Big)  } \ \dt
\]
hold
for any $0 \leq \tau_1 < \tau_2 \leq T$ and any test function $\varphi \in \DC([0,T] \times \Omega)$, $\varphi \in \DC([0,T] \times \Omega; R^3)$, respectively,
where the potential $\vc{V}$ is the unique solution of the Poisson equation
\bFormula{i11}
\Delta \vc{V}(t, \cdot) = \vr(t, \cdot) - \Ov{\vr} , \ \intO{ \vc{V}(t, \cdot) } = 0,\ t \in [0,T],
\ \mbox{with} \ \Ov{\vr} = \frac{1}{|\Omega|} \intO{ \vr_0 }.
\eF

\eD

%}

\bRemark{i1}

In view of (\ref{i9-}), the vacuum set where $\vr = 0$ is of zero Lebesgue measure, in particular, all terms in the integral identities (\ref{i9}),
(\ref{i10}) are well defined. This is in good agreement with the interpretation of $\vr$ as the density of a quantum fluid given by (\ref{i3+}),
(\ref{i3++}) as the nodal
zones of the Schr\" odinger equation are likely to be composed of ``tiny'' sets, see Kenig et al. \cite{KePoVe}, Seo \cite{Seo}.

\eR

The present paper examines the well/ill posedness of the Euler-Korteweg-Poisson system in the class of weak solutions introduced above. Observe that
for the particular choice $\vr \equiv \Ov{\vr}$, the problem (\ref{i1}-\ref{i3}) reduces to the ``damped'' Euler system with \emph{zero} pressure. In view of the recent ground-breaking results by DeLellis and Sz\' ekelyhidi \cite{DelSze}, \cite{DelSze2}, \cite{DelSze3} based on the method of \emph{convex integration}, such a system is ill-posed in the class of
weak solutions, meaning it admits
infinitely many solutions for \emph{any} initial data.

Chiodaroli \cite{Chiod} obtained similar illposedness results for the compressible Euler system using a ``non-constant'' coefficient version of the method of
\cite{DelSze3}; later the method was further extended in \cite{ChiFeiKre} in order to attack the more complex Euler-Fourier system. The main idea, elaborated in \cite{ChiFeiKre},
is to consider the Helmholtz decomposition
\[
\vJ = \vc{v} + \Grad \Psi , \ \Div \vc{v} = 0,
\]
to determine $\vr$ along with the acoustic potential $\Psi$, and to ``solve'' the momentum equation for $\vJ$ as a ``pressureless'' Euler system with
nonconstant coefficients. Adapting this approach to the present problem features an essential difficulty related to the presence of vacuum zones, where
the equations become singular. To overcome this problem, we extend the technique of convex integration to problems with non-constant \emph{singular} coefficients. In particular, we show a variant of the crucial \emph{oscillatory increment} lemma on an arbitrary open set by means of a careful scale analysis of its original version in \cite{DelSze3} and an application of Whitney covering lemma.

The solutions obtained by the method of convex integration suffer the well-known deficit that eliminates most of them as physically irrelevant: Although  their energy remains bounded at any instant $t$ including $t = 0$, they do not satisfy the total energy balance (\ref{i8}), not even as an inequality. In particular, the energy at any positive time may become strictly larger than that of the initial data. This motivates introducing the energy inequality
\bFormula{i12}
\intO{ E(\vr, \Grad \vr, \vJ )(\tau, \cdot) } + \int_0^\tau \intO{ \frac{ |\vc{J}|^2 }{\vr} } \ \dt \leq
\intO{ E(\vr_0, \Grad \vr_0, \vJ_0 ) } \ \mbox{for a.a.}\ \tau \in (0,T)
\eF
as a suitable \emph{admissibility criterion} in the class of weak solutions. Indeed we show that the \emph{dissipative} weak solutions, meaning the
weak solutions satisfying (\ref{i12}), enjoy the weak-strong uniqueness property - they coincide with the strong solution emanating from the same initial data as long as the latter exists. This result will be a direct consequence of the method of relative entropies adapted from \cite{FeNoJi}, \cite{FeNoSun}.

Finally, we note that even the dissipative weak solution may fail to be unique, at least for certain (non-smooth) initial data. Such a result follows from
a refined application of convex integration in the spirit of DeLellis and Sz\' ekelyhidi \cite{DelSze3}.

The paper consists of two parts. In the first one, we discuss the problem of well/ill posedness of the Euler-Korteweg-Poisson system in the class of
weak solutions. We start by stating the main result on the existence of infinitely many solutions in Section \ref{m}. In Section \ref{c}, we show how the method of convex integration can be adapted to the present setting and reduce the problem to \emph{oscillatory lemma} proved in Section \ref{o}.
The second part concerns the dissipative weak solutions introduced in Section \ref{d}. In Section \ref{r}, we show that the dissipative weak solutions
possess the weak-strong uniqueness property. Finally, we discuss the ill posedness of the Euler-Korteweg-Poisson system in the class of
dissipative weak solutions for particular initial data.

\section{Well/ill posedness in the class of weak solutions}
\label{m}

We start by introducing certain technical assumptions imposed on the structural properties of the functions $p = p(\vr)$, $\chi = \chi(\vr)$, specifically,
\bFormula{m1}
p \in C^1[0, \infty) \cap C^2(0, \infty), \ p(0) = 0, \ \chi \in C^2[0, \infty), \ \chi > 0 \ \mbox{in}\ (0, \infty).
\eF
Note that, in view of possible applications to the theory of quantum fluids, the pressure $p$ need not be monotone, not even positive.

The main result of the first part of the paper reads:

%\Cbox{Cgrey}{

\bTheorem{m1}
Under the hypotheses (\ref{m1}), suppose that the initial data satisfy
\bFormula{m2}
\vr_0 = r_0^2, \ r_0 \in C^2(\Omega), \ {\rm meas} \left\{ x \in \Omega \ \Big| \ r_0(x) = 0 \right\} = 0,
\eF
\bFormula{m3}
\vJ_0 = \vr_0 \vc{U}_0, \ \vc{U}_0 \in C^3( \Omega; R^3).
\eF

Then the initial value problem (\ref{i1}-\ref{i3}), (\ref{i4}) admits infinitely many weak solutions in $(0,T)$ in the sense specified in Definition
\ref{Di1}.

\eT

%}

\bRemark{m1}

It is easy to check that the initial data satisfying (\ref{m2}), (\ref{m3}) possess uniformly bounded energy $E(\vr_0, \Grad \vr_0, \vJ_0)$.
The hypothesis (\ref{m3}) could be relaxed, the present form asserts the existence of the \emph{initial velocity} $\vc{U}_0$.

\eR

The following two sections will be devoted to the proof of Theorem \ref{Tm1}. We first extend the density $\vr$ to the whole time interval $[0,T]$ and then
construct the desired weak solutions by the method of convex integration.

\section{Convex integration}
\label{c}

We start by extending the initial data $\vr_0$, $\vJ_0$ as a suitable solution $[\vr , \tilde{ \vJ }]$ to the equation of continuity on the whole time interval $[0,T]$. The function $\vr = \vr(t,\cdot)$ will be the unique solution of the transport equation
\bFormula{c1}
\partial_t \vr + \Div (\vr [\vc{U}_0 - \vc{Z}] ) = \partial_t \vr + [\vc{U}_0 - \vc{Z} ] \cdot \Grad \vr + \vr \Div \vc{U}_0 = 0, \ \vr(0, \cdot) = \vr_0,
\eF
where the spatially homogeneous vector function $\vc{Z} = \vc{Z}(t)$ is chosen in such a way that
\bFormula{c1+}
e^t \intO{ \vr [\vc{U}_0 - \vc{Z}] } = \intO{ \vr_0 \vc{U}_0 } \ \mbox{for all}\ t \in [0,T],
\eF
in particular $\vc{Z}(0) = 0$.

Indeed, for any given $\vc{Z} \in C([0,T]; R^3)$, the Cauchy problem (\ref{c1}) admits a unique solution $\vr$ and we may define a mapping
\[
\mathcal{T}: \vc{Z} \mapsto  \left( \intO{ \vr_0 } \right)^{-1} \left( \intO{ \vr \vc{U}_0 } - e^{-t} \intO{ \vr_0 \vc{U}_0 } \right).
\]
Clearly, the satisfaction of (\ref{c1+}) corresponds to finding a fixed point of the mapping $\mathcal{T}$. To this end, it is enough to observe that
the maximum of $\vr$ satisfying (\ref{c1}) is \emph{independent} of $\vc{Z}$, and
\[
\partial_t \mathcal{T}[Z] = \left( \intO{ \vr_0 } \right)^{-1} \left( \intO{ \partial_t \vr \vc{U}_0 } + e^{-t} \intO{ \vr_0 \vc{U}_0 } \right)
\]
\[
= \left( \intO{ \vr_0 } \right)^{-1} \left( \intO{ \vr \Grad \vc{U}_0 \cdot [\vc{U}_0 - \vc{Z}] } + e^{-t} \intO{ \vr_0 \vc{U}_0 } \right);
\]
whence the existence of a fixed point $\vc{Z}$ follows by a direct application of the Schauder  theorem in a bounded ball of
$C([0,T]; R^3)$.

Since $\vc{U}_0$ enjoys the regularity
(\ref{m3}), we deduce that

\begin{itemize}
\item $\vr(t, \cdot) \in C^2(\Omega)$ for any $t \in [0,T]$;
\item
\bFormula{c2}
{\rm meas} \left\{ x \in \Omega \ \Big| \ \vr(t, x) = 0 \right\} = 0 \ \mbox{for any}\ t \in [0,T];
\eF
\item for $\tilde {\vJ}(t,x) = \vr(t,x) \Big( \vc{U}_0(x) - \vc{Z}(t) \Big)$ we have
\bFormula{c3-}
e^t \intO{ \tilde{ \vJ }(t, \cdot) } = \intO{ \vJ_0 } \ \mbox{for any}\ t \in [0,T],
\eF
and
\bFormula{c3}
E(\vr, \Grad \vr, \tilde {\vc{J}} ) (t,\cdot) \leq \Ov{E} \ \mbox{for all}\ t \in [0,T].
\eF
\end{itemize}

\bRemark{c1}

Let $\vc{H}$ denote the standard Helmholtz projection onto the space of solenoidal functions. We have
\[
\intO{ \tilde{ \vJ } } = \intO{ \vc{H} [ \tilde{\vJ} ]},
\]
and (\ref{c3-}) yields
\bFormula{c3--}
\partial_t \intO{ \vc{H} [ \tilde{\vJ} ]} + \intO{ \vc{H} [ \tilde{\vJ} ]} = 0, \ \intO{ \vc{H} [ \tilde{\vJ} ](0, \cdot) } = \intO{ \vc{H} [ {\vJ}_0 ]}.
\eF
This relation is important in the construction of the so-called subsolutions introduced below.

\eR

\subsection{Convex integration ansatz}

The density $\vr$ being fixed through (\ref{c1}), we look for the flux $\vJ$ in the form
\[
\vJ = \vc{w} + \tilde{ \vJ },
\]
where
\bFormula{c4}
\vc{w} \in C_{\rm weak}([0,T], L^2(\Omega;R^3)) \cap L^\infty((0,T) \times \Omega;R^3), \ \Div \vc{w} = 0,\
\vc{w}(0, \cdot) = 0.
\eF
In particular, the equation of continuity (\ref{i9}), together with the initial conditions (\ref{i9+}), are satisfied.

In order to comply with (\ref{i10}), the function $\vc{w}$ must be taken such that
\bFormula{c5}
\partial_t \left( \vc{w} + \tilde{ \vJ } \right) + \Div \left( \frac{ ( \vc{w} + \tilde{ \vJ } ) \otimes (\vc{w} + \tilde {\vJ}) }{\vr} \right)
+ \Grad p(\vr) + \left( \vc{w} + \tilde{ \vJ } \right) =
\eF
\[
\Grad \Big( \chi (\vr) \Del \vr  \Big) + \frac{1}{2} \Grad \Big( \chi'(\vr) |\Grad \vr|^2 \Big) - 4 \Div \Big( \chi(\vr) \Grad \sqrt{\vr} \otimes \Grad \sqrt{\vr} \Big) + \vr \Grad V
\]
in the sense of distributions. Note that, in accordance with (\ref{c2}), (\ref{c3}), all quantities are bounded continuous functions on the (open) set
where $\vr > 0$, the complement of which in $\Omega$ is of zero measure.

For future analysis, it is convenient to rewrite (\ref{c5}) in a different form. We proceed in several steps:

\medskip

{\bf Step 1}

To begin, we write $\tilde{ \vJ }$ in terms of the Helmholtz projection as
\[
\tilde {\vJ} = \vc{H}[ \tilde{ \vJ } ] + \Grad M ;
\]
whence, replacing  $\vc{w} \approx \vc{w} + \vc{H}[ \tilde{ \vJ }]$, we convert (\ref{c4}), (\ref{c5}) to
\bFormula{c6}
\vc{w} \in C_{\rm weak}([0,T], L^2(\Omega;R^3)) \cap L^\infty((0,T) \times \Omega;R^3), \ \Div \vc{w} = 0,\
\vc{w}(0, \cdot) =  \vc{H}[\vJ_0],
\eF
\bFormula{c7}
\partial_t  \vc{w}  + \Div \left( \frac{ ( \vc{w} + \Grad M ) \otimes (\vc{w} + \Grad M) }{\vr} \right)
+  \vc{w} + \Grad \left( p(\vr) + \partial_t M  + M \right) =
\eF
\[
\Grad \Big( \chi (\vr) \Del \vr  \Big) + \frac{1}{2} \Grad \Big( \chi'(\vr) |\Grad \vr|^2 \Big) - 4 \Div \Big( \chi(\vr) \Grad \sqrt{\vr} \otimes \Grad \sqrt{\vr} \Big) + \vr \Grad V.
\]

\medskip

{\bf Step 2}

Multiplying (\ref{c7}) by $e^t$ and introducing a new quantity $\vc{v} = e^t \vc{w}$ we obtain
\bFormula{c8}
\vc{v} \in C_{\rm weak}([0,T], L^2(\Omega;R^3)) \cap L^\infty((0,T) \times \Omega;R^3), \ \Div \vc{v} = 0,\
\vc{v}(0, \cdot) = \vc{H}[\vJ_0],
\eF
\bFormula{c9}
\partial_t  \vc{v}  + \Div \left( \frac{ ( \vc{v} + e^t \Grad M ) \otimes (\vc{v} + e^t \Grad M) }{ e^t \vr} \right)
 + \Grad \left( e^t p(\vr) + e^t \partial_t M  + e^t M \right) =
\eF
\[
\Grad \Big( e^t \chi (\vr) \Del \vr  \Big) + \frac{1}{2} \Grad \Big( e^t \chi'(\vr) |\Grad \vr|^2 \Big) - 4 \Div \Big( e^t \chi(\vr) \Grad \sqrt{\vr} \otimes \Grad \sqrt{\vr} \Big) + e^t \vr \Grad V,
\]

\medskip

{\bf Step 3}

Finally, writing
\[
\vr \Grad V = \Div \left( \Grad V \otimes \Grad V - \frac{1}{3} |\Grad V|^2 \tn{I} \right) + \Grad \left( \Ov{\vr} V -\frac{1}{6} |\Grad V|^2 \right)
\]
and introducing new quantities
\[
r(t,x) = e^t \vr(t,x),
\]
\[
\vc{h}(t,x) = e^t \Grad M(t,x),
\]
\[
\Pi(t,x) = e^t \left( p(\vr) + \partial_t M + M - \chi(\vr) \Del \vr - \frac{1}{2} \chi'(\vr) |\Grad \vr|^2 + \frac{4}{3}
\chi(\vr) |\Grad \sqrt{\vr} |^2 - \Ov{\vr} V +\frac{1}{6} |\Grad V|^2 \right),
\]
\[
\tn{H}(t,x) = 4 e^t \left( \chi(\vr) \Grad \sqrt{\vr} \otimes \Grad \sqrt{\vr} - \frac{1}{3} \chi(\vr) |\Grad \sqrt{\vr} |^2 \tn{I} -
\frac{1}{4}\Grad V \otimes \Grad V + \frac{1}{12} |\Grad V|^2 \tn{I} \right),
\]
we obtain (\ref{c9}) in a concise form
\bFormula{c10}
\partial_t  \vc{v}  + \Div \left( \frac{ ( \vc{v} + \vc{h} ) \otimes (\vc{v} + \vc{h}) }{ r } + \tn{H} \right)
 + \Grad \Pi = 0,
\eF
where $\tn{H}$ is a symmetric traceless tensor.

\subsection{Subsolutions}

Let
\[
R_+ = \left\{ (t,x) \in (0,T) \times \Omega \ \Big| \ r(t,x) > 0 \right\}
\]
denote the set of positivity of the density $\vr$. In accordance with (\ref{c2}), $R_+$ is an open set of full measure in $(0,T) \times \Omega$.

Following DeLellis and Sz\' ekelyhidi \cite{DelSze3},   we introduce the set of subsolutions
\[
X_{0,e} = \left\{ \vc{v} \in C_{\rm weak}([0,T]; L^2(\Omega;R^3)) \ \Big| \ \vc{v}(0,\cdot) = \vc{H}[\vJ_0] \right.
\]
\[
\vc{v} \in C^1 ((0,T) \times \Omega;R^3), \ \partial_t \vc{v} + \Div \tn{U} = 0 \ \mbox{in}\ (0,T) \times \Omega
\ \mbox{for a certain}\ \tn{U} \in C^1((0,T) \times \Omega; R^{3 \times 3}_{\rm sym,0} ) ,
\]
\[
\left. \frac{3}{2} \lambda_{\rm max} \left[ \frac{( \vc{v} + \vc{h}) \otimes (\vc{v} + \vc{h} ) }{r} + \tn{H} - \tn{U} \right] < e
\ \mbox{in} \ R_+ \right\},
\]
where $\lambda_{\rm max}[\tn{A}]$ stands for the maximal eigenvalue of a symmetric matrix $\tn{A}$.

Here, the functions $\vc{h}$, $r$, $\tn{H}$ are the same as in (\ref{c10}), whereas the ``energy'' $e$ is taken in the form
\bFormula{c11-}
e(t,x) = \omega(t) - \frac{3}{2} \Pi(t,x),
\eF
where $\Pi$ is the ``pressure'' in (\ref{c10}) while $\omega$ is a suitable spatially homogeneous function specified below.
In accordance with (\ref{c3}), we have
\bFormula{c11}
\Pi \in L^\infty((0,T) \times \Omega), \ \Pi \in C(R_+).
\eF

Finally, seeing that
\[
\frac{3}{2} \lambda_{\rm max} \left[ \frac{( \vc{v} + \vc{h}) \otimes (\vc{v} + \vc{h} ) }{r} + \tn{H} - \tn{U} \right] \geq \frac{1}{2}
\frac{ | \vc{v} + \vc{h} |^2 }{r}
\]
we introduce a non-positive functional
\bFormula{c12}
\mathcal{I} [\vc{v}] = \int_{R_+} \left( \frac{1}{2} \frac{ | \vc{v} + \vc{h} |^2 }{r} - e \right) \ \dxdt =
\int_0^T \intO{ \left( \frac{1}{2} \frac{ | \vc{v} + \vc{h} |^2 }{r} - e \right) } \ \dt.
\eF

\subsection{Proof of Theorem \ref{Tm1} by convex integration}

The crucial ingredient of the proof of Theorem \ref{Tm1} is the following \emph{oscillatory lemma}.

%\Cbox{Cgrey}{

\bLemma{osc}

Let $U \subset R \times R^3$ be a bounded open set. Suppose that
\[
\vc{g} \in C(U; R^3), \ \tn{W} \in C(U; R^{3 \times 3}_{{\rm sym},0}), \ e,r \in C(U),\ r > 0, \ e \leq \Ov{e} \ \mbox{in}\ U
\]
are given such that
\[
\frac{3}{2} \lambda_{\rm max} \left[ \frac{ \vc{g} \otimes \vc{g} }{r} - \tn{W} \right] < e  \ \mbox{in}\
U.
\]

Then there exist sequences
\[
\vc{w}_n \in \DC (U;R^3), \ \tn{V}_n \in \DC(U; R^{3 \times 3}_{\rm sym,0}), \ n = 0,1,\dots
\]
such that
\[
\partial_t \vc{w}_n + \Div \tn{V}_n = 0 , \ \Div \vc{w}_n = 0 \ \mbox{in}\ R^3,
\]
\[
\frac{3}{2} \lambda_{\rm max} \left[ \frac{ (\vc{g} + \vc{w}_n) \otimes (\vc{g} + \vc{w}_n)}{r} - (\tn{W} + \tn{V}_n) \right] < e  \ \mbox{in}\ U,
\]
and
\[
\liminf_{n \to \infty} \int_{U} \frac{ | \vc{w}_n |^2 }{r} \ \dxdt \geq c(\Ov{e}) \int_{U} \left( e - \frac{1}{2} \frac{|\vc{g} |^2}{r} \right)^2 \ \dxdt.
\]

\eL

%}

\bRemark{m2}

We point out that the functions $\vc{g}$, $\tn{W}$ are continuous but not necessarily bounded on the open set $U$. Similarly, $r$ need not be bounded below away from zero. Thus Lemma \ref{Losc} can be interpreted as a singular version of similar results in \cite{Chiod}, \cite{DelSze3}.

\eR

The proof of Lemma \ref{Losc} will be given in Section \ref{o}. Taking this result for granted, the proof of Theorem \ref{Tm1} follows the
arguments similar to \cite{DelSze3}:

\medskip

{\bf Step 1}

First we observe that the set of subsolutions $X_{0,e}$ is non-empty, at least for a sufficiently large function $\omega$ in (\ref{c11-}). To see this, it is enough to take
\[
\vc{v} = e^t \vc{H}[ \tilde{ \vJ } ].
\]

As a consequence of (\ref{c3}) we get
\[
\frac{ |\vc{v} + \vc{h}|^2 }{r} = e^{t} \frac{ \left| \vc{H}[ \tilde{ \vJ } ] + \Grad M \right|^2 }{\vr} =
e^t   \frac{ |\tilde {\vJ} |^2 }{\vr} \leq \Ov{E}.
\]

Thus it is enough to find a suitable field $\tn{U} \in C^1((0,T) \times \Omega; R^{3 \times 3}_{{\rm sym},0})$ such that
\[
\Div \tn{U} = - \partial_t \vc{v}.
\]
This can be achieved by solving, for instance, the elliptic system
\[
\Div \left( \Grad \vc{w} + \Grad^t \vc{w} - \frac{2}{3} \Div \vc{w} \tn{I} \right) = - \partial_t \vc{v}
\]
since we have, by virtue of (\ref{c3--}),
\[
\intO{ \partial_t \vc{v} } = \partial_t \intO{ e^t \vc{H}[\tilde{ \vJ }] } = 0.
\]

As $\Pi (t,x)$ is bounded, we can choose $\omega$ in (\ref{c11-}) so large that $\vc{v} \in X_{0,e}$.

\medskip

{\bf Step 2}

Applying \emph{oscillatory lemma} (Lemma \ref{Losc}) with
\[
U = R_+, \ \vc{g} = \vc{h} + \vc{v}, \ \tn{W} = \tn{H} - \tn{U}
\]
we deduce that cardinality of the space $X_{0,e}$ is infinite.

\medskip

{\bf Step 3}

The last step leans on a sophisticated Baire category argument due to DeLellis and Sz\' ekelyhidi \cite{DelSze3}.
Endowing $X_{0,e}$ with the metrizable (on $X_{0,e}$) topology of the space $C_{\rm weak}([0,T]; L^2(\Omega;R^3))$, we deduce, by means of
that the functional
$\mathcal{I}$, defined through (\ref{c12}), admits infinitely many points of continuity on the closure of $X_{0,e}$ satisfying
\bFormula{ome}
\omega - \frac{3}{2} \Pi = e = \frac{1}{2} \frac{ |\vc{v} + \vc{h} |^2 }{r} , \
\tn{U} = \tn{H} + \frac{ (\vc{v} + \vc{h} ) \otimes (\vc{v} + \vc{h} )}{r} - \frac{1}{3} \frac{ |\vc{v} + \vc{h} |^2 }{r},
\eF
\[
\partial_t \vc{v} + \Div \tn{U} = 0 \ \mbox{in the distributional sense},\ \vc{v}(0, \cdot) = \vc{H}[ \vJ_0 ],
\]
which is exactly equation (\ref{c10}), cf. \cite{DelSze3}.

We have proved Theorem \ref{Tm1}.

\section{Oscillatory lemma}
\label{o}

Our goal in this section is to prove Lemma \ref{Losc}.

\subsection{Basic result}

\label{BR}

We start with the following basic result due to DeLellis and Sz\' ekelyhidi (cf. also Chiodaroli \cite{Chiod}).

\bLemma{DS}
Let
\[
Q = \left\{ (t,x) \ \Big| \ t \in (0,1), \ x \in (0,1)^3 \right\}, \
\vc{v} \in R^3,\ \tn{U} \in R^{3 \times 3}_{{\rm sym},0} , \ e > 0
\]
satisfying
\bFormula{HYP1}
\frac{3}{2} \lambda_{\rm max} \left[ \vc{v} \otimes \vc{v} - \tn{U} \right] < e \leq \overline{e}.
\eF

Then there exists sequences
\[
\vc{w}_n \in \DC (Q;R^3), \ \tn{V}_n \in \DC(Q; R^{3 \times 3}_{\rm sym,0}), \ n = 0,1,\dots
\]
such that
\[
\partial_t \vc{w}_n + \Div \tn{V}_n = 0,  \ \Div \vc{w}_n = 0 \ \mbox{in}\ R^3,
\]
\bFormula{CON1}
\frac{3}{2} \lambda_{\rm max} \left[ (\vc{v} + \vc{w}_n) \otimes (\vc{v} + \vc{w}_n) - (\tn{U} + \tn{V}_n) \right] < e \leq \overline{e} \ \mbox{in}\ Q.
\eF
\[
\vc{w}_n \to 0 \ \mbox{weakly in} \ L^2 (Q; R^3),
\]
\bFormula{C5}
\liminf_{n \to \infty} \int_Q | \vc{w}_n |^2 \ \dxdt \geq c(\Ov{e}) \int_Q \left( e - \frac{1}{2} |\vc{v} |^2 \right)^2 \ \dxdt
\eF

\eL

\bRemark{o1}

It is important that the constant in (\ref{C5}) is independent of $e$, $\vc{v}$, and $\tn{U}$.

\eR

\subsection{Extending by scaling}

Rescaling $\vc{w}_n \approx \vc{w_n}(t/L, x/L)$, $\tn{V}_n = \tn{V}_n(t/L,x/L)$ we can extend the validity of Lemma \ref{Losc} to an arbitrary cube
\[
Q_L = LQ = \left\{ (t,x) \ \Big| \ t \in (0,L), \ x \in (0,L)^3 \right\}, \ L > 0,
\]
with the same constant $c(\Ov{e})$ in (\ref{C5}).

Now, using additivity of the integral, we observe that Lemma \ref{Losc} holds on any domain
\[
Q_{T,L} = \left\{ (t,x) \ \Big| \ t \in (0,T), \ x \in (0,L)^3 \right\}, \ T,L > 0,
\]
via a decomposition of $Q_{T,L}$ on a (finite) number of cubes.

Finally, introducing a new scaling
\[
\vc{w}_n \approx \sqrt{r}  {\vc{w}} ( t/\sqrt{r}, x), \ \tn{V}_n \approx \sqrt{r}  {\vc{V}}_n (t/\sqrt{r}, x), \ t \in (0,T), \ x \in (0,L)^3,
\]
for a positive constant $r$,
we conclude that the hypothesis (\ref{HYP1}) may be replaced by
\[
\frac{3}{2} \lambda_{\rm max} \left[ \frac{\vc{v} \otimes \vc{v}}{r} - \tn{U} \right] < e \leq \overline{e},
\]
with the conclusion of Lemma \ref{Losc} valid with the obvious changes
\bFormula{sss1}
\frac{3}{2} \lambda_{\rm max} \left[ \frac{ (\vc{v} + \vc{w}_n) \otimes (\vc{v} + \vc{w}_n)}{r} - (\tn{U} + \tn{V}_n) \right] < e \leq \overline{e},
\eF
\[
\liminf_{n \to \infty} \int_{Q_{T,L}} \frac{ | \vc{w}_n |^2 }{r} \ \dxdt \geq c(\Ov{e}) \int_{Q_{T,L}} \left( e - \frac{1}{2} \frac{|\vc{v} |^2}{r} \right)^2 \ \dxdt
\]
in (\ref{CON1}), (\ref{C5}), respectively.

\subsection{Continuous perturbation}

\label{CP}

Our goal is to extend Lemma \ref{Losc} to the case, where $\vc{v}$, $\tn{U}$, $r > 0$, and $e$ are \emph{continuous} functions on the (closed) cube
$\Ov{Q}_L$ satisfying
\[
\frac{3}{2} \lambda_{\rm max} \left[ \frac{\vc{v} \otimes \vc{v}}{r} - \tn{U} \right] < e \leq \overline{e} \ \mbox{in} \ \Ov{Q}_L.
\]
Let us point out that $r$, being continuous on $\Ov{Q}_L$, is bounded below away from zero.

Now, choosing $\delta > 0$ small enough we
decompose
\[
\Ov{Q}_L = \cup_{i=1}^m \Ov{Q}^i, \ Q^i \cap Q^j = \emptyset \ \mbox{for}\ i \ne j,
\]
where $Q^i$ are cubes that can be taken small enough so that
\bFormula{D1}
\frac{3}{2} \lambda_{\rm max} \left[ \frac{\vc{v}_i \otimes \vc{v}_i}{r_i} - \tn{U}_i \right] < e_i - \delta \ \mbox{in}\ Q^i,\ i = 1 , \dots, m
\eF
for arbitrary constant quantities
\[
\vc{v}_i = \vc{v}(t_{i,v},x_{i,v}), \ r_i = r (t_{i,r}, x_{i,r}), \tn{U}_i = \tn{U} (t_{i,u}, r_{i,u}), \ e_i = e (t_{i,r}, x_{i,r}),
\
(t_{i,\cdot},x_{i,\cdot}) \in Q^i.
\]
Moreover, by the same token, we may assume that
\bFormula{D2}
\left| \frac{3}{2} \lambda_{\rm max} \left[ \frac{ (\vc{v}_1 + \vc{w}) \otimes (\vc{v}_1 + \vc{w})}{r_1} - (\tn{U}_1 + \tn{V}) \right]
- \frac{3}{2} \lambda_{\rm max} \left[ \frac{ (\vc{v}_2 + \vc{w}) \otimes (\vc{v}_2 + \vc{w})}{r_2} - (\tn{U}_2 + \tn{V}) \right]
\right| < \frac{\delta}{2}
\eF
provided $\vc{w}$, $\tn{V}$ are bounded and
\[
\vc{v}_j = \vc{v}(t_{j,v},x_{j,v}), \ r_j = r (t_{j,r}, x_{j,r}), \tn{U}_j = \tn{U} (t_{j,u}, r_{j,u}),
\
(t_{j,\cdot},x_{j,\cdot}),  \in Q^i, \ j = 1,2.
\]

Thus, using (\ref{D1}), (\ref{D2}), together with the result for the constant coefficients shown above, we obtain the desired sequences
$\{ \vc{w}_n \}_{n=1}^\infty$, $\{ \tn{V}_n \}_{n=1}^\infty$ satisfying
\[
\liminf_{n \to \infty} \int_{Q_{T,L}} \frac{ | \vc{w}_n |^2 }{r} \ \dxdt \geq c(\Ov{e}) \int_{Q_{T,L}} \left( e - \delta - \frac{1}{2} \frac{|\vc{v} |^2}{r} \right)^2 \ \dxdt.
\]
As $\delta > 0$ can be taken arbitrarily small, we conclude.

\subsection{A decomposition lemma and the final result}

To conclude, we make use of the standard Whitney decomposition lemma, see Stein \cite{STEIN1}:

\bLemma{WH}

Let $U \subset R^N$ be an arbitrary open set. The there exists a countable family of (dyadic) open cubes $Q^i$ such that
\[
U = \cup_{i=1}^\infty \Ov{Q}^i, \ Q^i \cap Q^j = \emptyset \ \mbox{for}\ i \ne j,
\]
and
\bFormula{dist}
{\rm diam}[ Q^i ] \leq {\rm dist} [ Q^i, \partial U] \leq 4 {\rm diam} [Q^i] \ \mbox{for all}\ i = 1,,\dots
\eF

\eL

Decomposing the domain $U$ in Lemma \ref{Losc} as in Lemma \ref{LWH} and using the results of Section \ref{CP} on each cube $Q^i$, we complete the proof of Lemma \ref{Losc}. Note that, thanks to (\ref{dist}), the restriction of the continuous function $r$ to $Q^i$ is bounded below away from zero.

\section{Dissipative weak solutions}
\label{d}

A weak solution $\vr$, $\vJ$ is called \emph{dissipative weak solution} if, in addition to the stipulations listed in Definition \ref{Di1},
it satisfies the energy inequality (\ref{i12}). It is worth revisiting the weak solutions
constructed in the proof of Theorem \ref{Tm1} in the light of (\ref{i12}). We recall that the energy of these reads
\[
E(\vr, \Grad \vr, \vJ) = \frac{1}{2} \frac{ |\vc{J}|^2 }{\vr} + P(\vr) + 2 \chi(\vr)| \Grad \sqrt{\vr} |^2  + \frac{1}{2} |\Grad V|^2,
\]
where, by virtue of (\ref{c11-}), (\ref{ome}),
\[
\frac{1}{2} \frac{ |\vc{J}|^2 }{\vr} = e^{-t} \omega (t) - \frac{3}{2} \Pi (t,x).
\]
The function $\omega(t)$ could be chosen arbitrary but large enough, here large means in terms of the initial data. Going back to the energy inequality (\ref{i12})
we get
\[
\frac{{\rm d}}{{\rm d}t} \intO{ E(\vr, \Grad \vr, \vJ) } + \intO{ \frac{ |J| ^2 }{\vr} } = e^{-t} |\Omega| \Big( \omega'(t) + \omega(t) \Big) + h(t)
\]
for a certain function $h$ depending solely on the initial data. Thus, taking
\[
\omega (t) = e^{-2t}M , \ M = M(T) \ \mbox{large enough},
\]
we obtain
\bFormula{ei}
\frac{{\rm d}}{{\rm d}t} \intO{ E(\vr, \Grad \vr, \vJ) } + \intO{ \frac{ |J| ^2 }{\vr} } = - e^{-3t} |\Omega| M + h(t) \leq 0 \ \mbox{in} \ (0,T).
\eF
We conclude that the weak solutions in Theorem \ref{Tm1} can be constructed to satisfy the energy inequality in the \emph{open} interval. On the other hand,
in general, they are not expected to satisfy (\ref{i12}), meaning the energy balance may be violated at the initial time $t = 0$. We come back to this issue at the end of this section.

\subsection{Relative energy (entropy) inequality}

To simplify the forthcoming presentation, we restrict ourselves to the case of constant ``capillarity'' tensor $K = 1$ or, equivalently,  $\chi(\vr) = {\vr}$.
Motivated by the analysis in \cite{FeNoJi}, we introduce \emph{relative energy} functional
\bFormula{d1}
\mathcal{E} \left( \vr,  \vJ \ \Big| \ r,  \vc{L} \right)
\eF
\[
=
\intO{ \left[ \frac{1}{2} \vr \left| \frac{\vJ}{\vr} - \frac{\vc{L} }{r} \right|^2   + P(\vr) - P'(r) (\vr - r) - P(r)   +
\frac{1}{2} | \Grad \vr - \Grad r |^2  \right] },
\]
where $r > 0$, $\vL$ are smooth functions.
Note that the expression
\[
\frac{1}{2} \vr \left| \frac{\vJ}{\vr} - \frac{\vc{L} }{r} \right|^2 = \frac{1}{2} \frac{|\vJ|^2 }{\vr} - {\vJ} \cdot \frac{\vc{L}}{r} +
\frac{1}{2} \vr \frac{ |\vc{L} |^2 }{r^2}
\]
makes sense for any finite energy weak solution to the Euler-Korteweg-Poisson system.

Similarly to \cite{FeNoJi}, \cite{FeNoSun}, we derive an inequality describing the time evolution of $\mathcal{E}$.

\medskip

{\bf Step 1}

Taking $\vL/r$ as a test function in the momentum balance (\ref{i10}) we obtain
\bFormula{d2}
- \intO{ \vJ \cdot \frac{\vL}{r} } \Big|^{t = \tau_2}_{t = \tau_1} =
- \int_{\tau_1}^{\tau_2} \intO{ \left[ \vJ \cdot \partial_t \left( \frac{\vL}{r} \right) + \frac{ \vJ \otimes \vJ }{\vr} : \Grad \left( \frac{\vL}{r} \right) + p(\vr) \Div \left( \frac{\vL}{r} \right) \right] } \ \dt
\eF
\[
- \int_{\tau_1}^{\tau_2} \intO{ \left[ \Grad \vr \cdot \Grad \left[ {\vr}  \Div \left( \frac{\vL}{r} \right) \right]
- \frac{1}{2} |\Grad \vr|^2 \Div \left( \frac{\vL}{r} \right) +\Grad {\vr} \otimes \Grad {\vr} : \Grad \left( \frac{\vL}{r} \right) \right] }
\]
\[
+ \int_{\tau_1}^{\tau_2} \intO{ \left[ \vc{J} \cdot \left( \frac{\vL}{r} \right) - \vr \Grad \vc{V} \cdot \left( \frac{\vL}{r} \right) \right]  } \ \dt.
\]

\medskip

{\bf Step 2}

Similarly, the choice $\varphi = \frac{1}{2} |\vL|^2/r$ in (\ref{i9}) gives rise to
\bFormula{d3}
\frac{1}{2} \intO{ \vr \frac{ |\vL|^2 }{r^2} } \Big|^{t = \tau_2}_{t = \tau_1} =
\frac{1}{2} \int_{\tau_1}^{\tau_2} \intO{ \left[ \vr \partial_t \left( \frac{ |\vL|^2 }{r^2}   \right) + \vJ \cdot \Grad \left( \frac{ |\vL|^2 }{r^2} \right) \right] } \ \dt.
\eF

\medskip

{\bf Step 3}

Taking $\varphi = P'(r)$ in (\ref{i9}) we get
\bFormula{d4}
- \intO{ \vr P'(r) } \Big|^{t = \tau_2}_{t = \tau_1} =
- \int_{\tau_1}^{\tau_2} \intO{ \left[ \vr \partial_t P'(r)  + \vJ \cdot \Grad P'(r) \right] } \ \dt.
\eF

\medskip

{\bf Step 4}

Finally, the test function $\varphi =  \Del r$ in (\ref{i9}) yields
\bFormula{d4a}
 \intO{ \Grad \vr \cdot \Grad r  } \Big|^{t = \tau_2}_{t = \tau_1} +
\int_{\tau_1}^{\tau_2} \intO{ \left( \vr \partial_t \Del r + \vJ \cdot \Grad \Del r \right) } \ \dt = 0,
\eF

\medskip

{\bf Step 5}

Summing up (\ref{ei}), (\ref{d2} - \ref{d4}) we obtain the \emph{relative energy inequality}
\[
\mathcal{E} \left( \vr, \vJ \ \Big| \ r , \vL \right) \Big|^{t = \tau_2}_{t = \tau_1} = \intO{ \left[ \frac{1}{2} \frac{ |\vJ|^2 }{\vr} - {\vJ} \cdot \frac{\vL}{r} + \frac{1}{2} \vr \frac{| \vL|^2}{r} + P(\vr) - P'(r) \vr + p(r) + \frac{1}{2} |\Grad \vr - \Grad r |^2 \right] } \Big|^{t = \tau_2}_{t = \tau_1}
\]
\[
\leq - \frac{1}{2} \intO{  |\Grad V |^2  } \Big|^{t = \tau_2}_{t = \tau_1} - \int_{\tau_1}^{\tau_2}
\intO{ \frac{ |\vc{J}|^2 }{\vr} } \ \dt + \int_{\tau_1}^{\tau_2} \intO{ \left( \vr \partial_t \Del r + \vJ \cdot \Grad \Del r - \partial_t r \Del r \right) } \ \dt
\]
\[
- \int_{\tau_1}^{\tau_2} \intO{ \left[ \vJ \cdot \partial_t \left( \frac{\vL}{r} \right) + \frac{ \vJ \otimes \vJ }{\vr} : \Grad \left( \frac{\vL}{r} \right) + p(\vr) \Div \left( \frac{\vL}{r} \right) \right] } \ \dt
\]
\[
- \int_{\tau_1}^{\tau_2} \intO{ \left[   \Grad \vr \cdot \Grad \left[ {\vr}  \Div \left( \frac{\vL}{r} \right) \right]
- \frac{1}{2} |\Grad \vr|^2 \Div \left( \frac{\vL}{r} \right) +\Grad {\vr} \otimes \Grad {\vr} : \Grad \left( \frac{\vL}{r} \right) \right] }
\]
\[
+ \int_{\tau_1}^{\tau_2} \intO{ \left[ \vc{J} \cdot \left( \frac{\vL}{r} \right) - \vr \Grad \vc{V} \cdot \left( \frac{\vL}{r} \right) \right]  } \ \dt
+
\frac{1}{2} \int_{\tau_1}^{\tau_2} \intO{ \left[ \vr \partial_t \left( \frac{ |\vL|^2 }{r^2}   \right) + \vJ \cdot \Grad \left( \frac{ |\vL|^2 }{r^2} \right) \right] } \ \dt
\]
\[
- \int_{\tau_1}^{\tau_2} \intO{ \left[ \vr \partial_t P'(r)  + \vJ \cdot \Grad P'(r) - \partial_t p(r) \right] } \ \dt
\]
that can be written in a more concise form:

\bFormula{d5}
\mathcal{E} \left( \vr, \vJ \ \Big| \ r , \vL \right) \Big|^{t = \tau_2}_{t = \tau_1} +
\int_{\tau_1}^{\tau_2} \intO{ \left( \vJ - \vr \frac{\vL}{r} \right) \cdot \left( \frac{\vJ}{\vr} - \frac{\vL}{r} \right) } \ \dt
\eF
\[
\leq - \frac{1}{2} \intO{  |\Grad V |^2  } \Big|^{t = \tau_2}_{t = \tau_1} - \int_{\tau_1}^{\tau_2}
\intO{ \vr \frac{\vL}{r} \cdot \left( \frac{\vJ}{\vr} - \frac{\vL}{r} \right) } \ \dt
\]
\[
+ \int_{\tau_1}^{\tau_2} \intO{ \left( \vr \partial_t \Del r + \vJ \cdot \Grad \Del r - \partial_t r \Del r \right) } \ \dt
- \int_{\tau_1}^{\tau_2} \intO{ \vr \Grad {V} \cdot \left( \frac{\vL}{r} \right)  } \ \dt
\]
\[
+ \int_{\tau_1}^{\tau_2} \intO{ \left[  \partial_t \left( \frac{\vL}{r} \right) + \frac{ \vJ }{\vr} \cdot \Grad \left( \frac{\vL}{r} \right) \right] \cdot
\left( \frac{\vr}{r} \vL - \vJ \right)   } \ \dt
\]
\[
- \int_{\tau_1}^{\tau_2} \intO{ \left[   \frac{1}{4} \Grad \vr \cdot \Grad \left( {\vr}  \Div \left( \frac{\vL}{r} \right) \right)
- \frac{1}{2} |\Grad \vr|^2 \Div \left( \frac{\vL}{r} \right) +\Grad {\vr} \otimes \Grad {\vr} : \Grad \left( \frac{\vL}{r} \right) \right] }
\]
\[
+ \int_{\tau_1}^{\tau_2} \intO{ \left[ (r - \vr) \partial_t P'(r)  + (\vL - \vJ) \cdot \Grad P'(r) + (p(r) - p(\vr)) \Div \left( \frac{\vL}{r} \right) \right] } \ \dt.
\]

\bRemark{lions}

The relative energy inequality (\ref{d5}) holds for any smooth ``test'' functions $r$, $\vc{L}$, $r$ bounded bellow away from zero. Alternatively, we may use (\ref{d5}) as a \emph{definition} of a dissipative solution in the spirit of Lions \cite{LI} in the context of the incompressible Euler system.

\eR

\subsection{Weak-strong uniqueness}
\label{r}

Supposing that the Euler-Korteweg-Poisson sytem  admits a \emph{smooth} solution $\tilde \vr > 0$, $\tilde{ \vJ }$, we set
$r = \tilde \vr$, $\vL = \tilde{ \vJ}$ as test functions in the relative energy inequality (\ref{d5}) to obtain
\bFormula{r1}
\mathcal{E} \left( \vr, \vJ \ \Big| \ \tilde \vr , \tilde{ \vJ } \right) \Big|^{t = \tau_2}_{t = \tau_1}
\leq - \frac{1}{2} \intO{  |\Grad V |^2  } \Big|^{t = \tau_2}_{t = \tau_1}
\eF
\[
+ \int_{\tau_1}^{\tau_2} \intO{ \left( \vr \partial_t \Del \tilde \vr + \vJ \cdot \Grad \Del \tilde \vr - \partial_t \tilde \vr \Del \tilde \vr \right) } \ \dt
\]
\[
+ \int_{\tau_1}^{\tau_2} \intO{ \vr
\left( \frac{ \vJ }{\vr} - \frac{ \tilde {\vJ} }{\tilde \vr} \right)  \cdot \Grad \left( \frac{\tilde{ \vJ }}{\tilde \vr} \right)  \cdot
\left( \frac{\tilde {\vJ} }{\tilde \vr}  - \frac{\vJ}{\vr} \right)   } \ \dt + \int_{\tau_1}^{\tau_2} \intO{
\Big( p(\tilde \vr ) - p'(\tilde \vr) (\tilde \vr - \vr) - p(\vr) \Big) \Div \left( \frac{\tilde{\vJ}}{\tilde \vr} \right) } \ \dt
\]
\[
+ \int_{\tau_1}^{\tau_2} \intO{ \left( \frac{\vr}{\tilde \vr} \tilde {\vJ} - \vJ \right) \cdot
\Big( \Grad \Del \tilde \vr +  \Grad \Del^{-1} (\tilde \vr - \Ov{\vr} ) \Big) } \ \dt - \int_{\tau_1}^{\tau_2} \intO{ \vr \Grad {V} \cdot \left( \frac{\tilde {\vJ} }{\tilde \vr} \right)  } \ \dt
\]
\[
- \int_{\tau_1}^{\tau_2} \intO{ \left[  \Grad \vr \cdot \Grad \left( {\vr} \Div \left( \frac{\tilde {\vJ }}{\tilde \vr} \right) \right)
  -\frac{1}{2} |\Grad \vr|^2 \Div \left( \frac{\tilde{ \vJ}}{\tilde \vr} \right) + \Grad {\vr} \otimes \Grad {\vr} : \Grad \left( \frac{\tilde {\vJ} }{\tilde \vr} \right) \right] }
\]
\[
+ \int_{\tau_1}^{\tau_2} \intO{ \left[ (\tilde \vr - \vr) \partial_t P'(\tilde \vr)  + (\tilde \vr - \vr) \frac{\tilde {\vJ}}{\tilde \vr} \cdot \Grad P'(\tilde \vr) + p'(\tilde \vr )(\tilde \vr - \vr)  \Div \left( \frac{\tilde{\vJ}}{\tilde \vr} \right) \right] } \ \dt,
\]
where, furthermore,
\[
(\tilde \vr - \vr) \partial_t P'(\tilde \vr)  + (\tilde \vr - \vr) \frac{\tilde {\vJ}}{\tilde \vr} \cdot \Grad P'(\tilde \vr) + p'(\tilde \vr )(\tilde \vr - \vr)  \Div \left( \frac{\tilde{\vJ}}{\tilde \vr} \right)
=
(\tilde \vr - \vr) P''(\tilde \vr) \Big( \partial_t \tilde \vr + \Div \tilde {\vJ } \Big) = 0,
\]
and
\[
\int_{\tau_1}^{\tau_2} \intO{ \left[ \left( \frac{\vr}{\tilde \vr} \tilde {\vJ} - \vJ \right) \cdot \Grad \Del^{-1} (\tilde \vr - \Ov{\vr} )
- \frac{\vr}{\tilde \vr} \tilde {\vJ} \cdot \Grad \Del^{-1} (\vr - \Ov{\vr}) \right] } \ \dt
\]
\[
= \int_{\tau_1}^{\tau_2} \intO{ \frac{\vr}{\tilde \vr} \tilde{ \vJ } \cdot \Grad \Del^{-1}( \tilde \vr - \vr ) } \ \dt
- \intO{ \vr \Del^{-1}( \tilde \vr - \Ov{\vr} ) } \Big|_{t = \tau_1}^{t = \tau_2} -
\int_{\tau_1}^{\tau_2} \intO{ \Del^{-1} (\vr - \Ov{\vr} ) \Div \tilde{ \vJ } } \ \dt
\]
\[
= \int_{\tau_1}^{\tau_2} \intO{ \frac{1}{\tilde \vr} \left( {\vr} - \tilde \vr \right) \tilde{ \vJ } \cdot \Grad \Del^{-1}( \tilde \vr - \vr ) } \ \dt
- \intO{ \vr \Del^{-1}( \tilde \vr - \Ov{\vr} ) } \Big|_{t = \tau_1}^{t = \tau_2} + \int_{\tau_1}^{\tau_2} \intO{
\partial_t \tilde \vr \Del^{-1} (\tilde \vr - \Ov{\vr} ) } \ \dt.
\]

Consequently, after a simple manipulation, the relation (\ref{r1}) reads
\bFormula{r2}
\mathcal{E} \left( \vr, \vJ \ \Big| \ \tilde \vr , \tilde{ \vJ } \right) \Big|^{t = \tau_2}_{t = \tau_1}
\leq \frac{1}{2} \intO{  (\tilde \vr - \vr) \Del^{-1} (\tilde \vr - \vr)  } \Big|^{t = \tau_2}_{t = \tau_1}
\eF
\[
+ \int_{\tau_1}^{\tau_2} \intO{ \vr
\left( \frac{ \vJ }{\vr} - \frac{ \tilde {\vJ} }{\tilde \vr} \right)  \cdot \Grad \left( \frac{\tilde{ \vJ }}{\tilde \vr} \right)  \cdot
\left( \frac{\tilde {\vJ} }{\tilde \vr}  - \frac{\vJ}{\vr} \right)   } \ \dt + \int_{\tau_1}^{\tau_2} \intO{
\Big( p(\tilde \vr ) - p'(\tilde \vr) (\tilde \vr - \vr) - p(\vr) \Big) \Div \left( \frac{\tilde{\vJ}}{\tilde \vr} \right) } \ \dt
\]
\[
+\int_{\tau_1}^{\tau_2} \intO{ \frac{1}{\tilde \vr} \left( {\vr} - \tilde \vr \right) \tilde{ \vJ } \cdot \Grad \Del^{-1}( \tilde \vr - \vr ) } \ \dt
\]
\[
+ \int_{\tau_1}^{\tau_2} \intO{  \frac{\vr}{\tilde \vr} \tilde {\vJ}  \cdot
\Grad \Del \tilde \vr  } \ \dt
 + \int_{\tau_1}^{\tau_2} \intO{ \left( \vr \partial_t \Del \tilde \vr  - \partial_t \tilde \vr \Del \tilde \vr \right) } \ \dt
\]
\[
- \int_{\tau_1}^{\tau_2} \intO{ \left[ \Grad \vr \cdot \Grad \left( {\vr} \Div \left( \frac{\tilde {\vJ }}{\tilde \vr} \right) \right)
  -\frac{1}{2} |\Grad \vr|^2 \Div \left( \frac{\tilde{ \vJ}}{\tilde \vr} \right) + \Grad {\vr} \otimes \Grad {\vr} : \Grad \left( \frac{\tilde {\vJ} }{\tilde \vr} \right) \right] }.
\]

The next step is formal in the sense that it requires higher regularity of $\vr$ but the final relation can be justified by means of a density argument. We write
\[
- \int_{\tau_1}^{\tau_2} \intO{ \left[  \Grad \vr \cdot \Grad \left( {\vr} \Div \left( \frac{\tilde {\vJ }}{\tilde \vr} \right) \right)
  -\frac{1}{2} |\Grad \vr|^2 \Div \left( \frac{\tilde{ \vJ}}{\tilde \vr} \right) + \Grad {\vr} \otimes \Grad {\vr} : \Grad \left( \frac{\tilde {\vJ} }{\tilde \vr} \right) \right] }
\]
\[
=
- \int_{\tau_1}^{\tau_2} \intO{ \frac{\vr}{\tilde \vr} \tilde{ \vJ} \cdot \Grad \Del  \vr} \ \dt,
\]
and, consequently,
\[
\int_{\tau_1}^{\tau_2} \intO{ \frac{\vr}{\tilde \vr} \tilde{ \vJ} \cdot \Grad \Del  (\tilde \vr - \vr) } \ \dt
\]
\[
 =
\int_{\tau_1}^{\tau_2} \intO{ \tilde{ \vJ} \cdot \Grad \Del  (\tilde \vr - \vr) } \ \dt +
\int_{\tau_1}^{\tau_2} \intO{ \frac{\tilde {\vJ }}{\tilde \vr} (\vr - \tilde \vr) \cdot \Grad \Del (\vr - \tilde \vr) } \ \dt
\]
\[
=
\int_{\tau_1}^{\tau_2} \intO{ \partial_t \tilde \vr \Del  (\tilde \vr - \vr) } \ \dt+
\int_{\tau_1}^{\tau_2} \intO{ \frac{\tilde {\vJ }}{\tilde \vr} (\vr - \tilde \vr) \cdot \Grad \Del (\vr - \tilde \vr) } \ \dt.
\]
Thus, going back to (\ref{r2}) we obtain
\bFormula{r3}
\mathcal{E} \left( \vr, \vJ \ \Big| \ \tilde \vr , \tilde{ \vJ } \right) \Big|^{t = \tau_2}_{t = \tau_1}
- \frac{1}{2} \intO{  (\tilde \vr - \vr) \Del^{-1} (\tilde \vr - \vr)  } \Big|^{t = \tau_2}_{t = \tau_1}
\eF
\[
\leq
\int_{\tau_1}^{\tau_2} \intO{ \vr
\left( \frac{ \vJ }{\vr} - \frac{ \tilde {\vJ} }{\tilde \vr} \right)  \cdot \Grad \left( \frac{\tilde{ \vJ }}{\tilde \vr} \right)  \cdot
\left( \frac{\tilde {\vJ} }{\tilde \vr}  - \frac{\vJ}{\vr} \right)   } \ \dt + \int_{\tau_1}^{\tau_2} \intO{
\Big( p(\tilde \vr ) - p'(\tilde \vr) (\tilde \vr - \vr) - p(\vr) \Big) \Div \left( \frac{\tilde{\vJ}}{\tilde \vr} \right) } \ \dt
\]
\[
+\int_{\tau_1}^{\tau_2} \intO{ \frac{\tilde{ \vJ }}{\tilde \vr} \left( {\vr} - \tilde \vr \right)  \cdot \Grad \Del^{-1}( \tilde \vr - \vr ) } \ \dt
+ \int_{\tau_1}^{\tau_2} \intO{ \frac{\tilde {\vJ }}{\tilde \vr} (\vr - \tilde \vr) \cdot \Grad \Del (\vr - \tilde \vr) } \ \dt.
\]

Finally,
\[
\intO{ \frac{\tilde {\vJ }}{\tilde \vr} (\vr - \tilde \vr) \cdot \Grad \Del (\vr - \tilde \vr) } =
\]
\[
= \intO{ \frac{\tilde {\vJ }}{\tilde \vr} \cdot \Div \left[ (\vr - \tilde \vr)  \Del (\vr - \tilde \vr)  + \frac{1}{2}  |\Grad (\vr - \tilde \vr) |^2  \tn{I} -  \Grad (\vr - \tilde \vr) \otimes \Grad ( {\vr} - \tilde \vr ) \right] } \ \dt
\]
\[
= \intO{ \Grad (\vr - \tilde \vr) \otimes \Grad ( {\vr} - \tilde \vr ) : \Grad \left( \frac{\tilde {\vJ }}{\tilde \vr} \right) } -
\frac{3}{2} \intO{ \Div \left( \frac{\tilde {\vJ }}{\tilde \vr} \right) |\Grad (\vr - \tilde \vr) |^2 }
+ \intO{ (\vr - \tilde \vr) \Grad \Div \left( \frac{\tilde {\vJ }}{\tilde \vr} \right) \cdot \Grad (\vr - \tilde \vr)}.
\]
Thus the relation (\ref{r3}) takes its final form
\bFormula{r4}
\mathcal{E} \left( \vr, \vJ \ \Big| \ \tilde \vr , \tilde{ \vJ } \right) \Big|^{t = \tau_2}_{t = \tau_1}
- \frac{1}{2} \intO{  (\tilde \vr - \vr) \Del^{-1} (\tilde \vr - \vr)  } \Big|^{t = \tau_2}_{t = \tau_1}
\eF
\[
\leq
\int_{\tau_1}^{\tau_2} \intO{ \vr
\left( \frac{ \vJ }{\vr} - \frac{ \tilde {\vJ} }{\tilde \vr} \right)  \cdot \Grad \left( \frac{\tilde{ \vJ }}{\tilde \vr} \right)  \cdot
\left( \frac{\tilde {\vJ} }{\tilde \vr}  - \frac{\vJ}{\vr} \right)   } \ \dt + \int_{\tau_1}^{\tau_2} \intO{
\Big( p(\tilde \vr ) - p'(\tilde \vr) (\tilde \vr - \vr) - p(\vr) \Big) \Div \left( \frac{\tilde{\vJ}}{\tilde \vr} \right) } \ \dt
\]
\[
+\int_{\tau_1}^{\tau_2} \intO{ \frac{\tilde{ \vJ }}{\tilde \vr} \left( {\vr} - \tilde \vr \right)  \cdot \Grad \Del^{-1}( \tilde \vr - \vr ) } \ \dt
\]
\[
\intO{ \Grad (\vr - \tilde \vr) \otimes \Grad ( {\vr} - \tilde \vr ) : \Grad \left( \frac{\tilde {\vJ }}{\tilde \vr} \right) } +
\frac{3}{2} \intO{ \Div \left( \frac{\tilde {\vJ }}{\tilde \vr} \right) |\Grad (\vr - \tilde \vr) |^2 }
+ \intO{ (\vr - \tilde \vr) \Grad \Div \left( \frac{\tilde {\vJ }}{\tilde \vr} \right) \cdot \Grad (\vr - \tilde \vr)}.
\]

Applying Gronwall's lemma to (\ref{r4}) we deduce that
\[
\mathcal{E} \left( \vr, \vJ \ \Big| \ \tilde \vr , \tilde{ \vJ } \right)(t, \cdot) = 0, \ t \in [0,T]
\ \mbox{as soon as} \ \mathcal{E} \left( \vr, \vJ \ \Big| \ \tilde \vr , \tilde{ \vJ } \right)(0, \cdot) = 0.
\]
We have shown the following \emph{weak-strong uniqueness property} of the Euler-Korteweg-Poisson system

%\Cbox{Cgrey}{

\bTheorem{ws}
Let $K(\vr) = \Ov{K}$ be a positive constant.
Let $\vr$, $\vJ$ be a \emph{dissipative} weak solution to the Euler-Korteweg-Poisson system (\ref{i1} - \ref{i3}), (\ref{i4}), with
in $(0,T) \times \Omega$ such that
\[
\vr(t,x) \geq \underline{\vr} > 0 \ \mbox{for a.a.}\ (t,x) \in (0,T) \times \Omega.
\]
Suppose that the problem (\ref{i1} - \ref{i3}), (\ref{i4}) admits a classical (strong) solution $\tilde \vr$, $\tilde{ \vJ }$
in $(0,T) \times \Omega$ emanating from the same initial data as $\vr$, $\vJ$.

Then
\[
\vr \equiv \tilde \vr,\ \vJ \equiv \tilde \vJ.
\]

\eT

%}

\bRemark{tzav2}

Having finished this paper, we have learned that more general results of the same type were obtained by Giesselmann, Lattanzio and Tzavaras \cite{GieLatTza}.

\eR

\bRemark{marc}

The \emph{existence} of local-in-time regular solutions and the \emph{existence} of global-in-time regular solutions where the initial data are taken as small perturbations of subsonic steady states to the quantum hydrodynamics system was proved by Li and one of the authors in \cite{LiMarcati}.

\eR

\subsection{Concluding remarks}
\label{cr}

Summarizing the previous discussion we may infer that:

\begin{itemize}

\item The Euler-Korteweg-Poisson system admits \emph{infinitely many} weak solutions for \emph{any} sufficiently smooth initial data. The solutions are defined on an arbitrary time interval $(0,T)$, where they satisfy the energy inequality, with a possible exception of the initial time $t= 0$.

\item The dissipative weak solution satisfy the energy as well as the relative energy inequality in $(0,T)$. They coincide with the strong solution emanating from the same initial data as long as the latter exists. In other words, the strong solutions are unique in the class of weak solutions.

\end{itemize}

In the light of the above arguments, it may seem plausible to eliminate the majority of the ``strange'' weak solutions obtained by the method of convex integration by stipulating the energy (relative energy) inequality. Unfortunately, however, a nowadays straightforward modification of the method
of convex integration yields the following result that can be proved, given the oscillatory lemma \ref{Losc}, by the arguments specified in
\cite{ChiFeiKre}.

%\Cbox{Cgrey}{

\bTheorem{cc1}
Let $\vr_0$ be given in the class specified in Theorem \ref{Tm1}. Let $T > 0$ be an arbitrary positive time. 

Then there exists the initial distribution of the momentum
\[
\vJ_0 \in L^\infty(\Omega; R^3)
\]
such that the Euler-Korteweg-Poisson system (\ref{i1} - \ref{i3}) admits infinitely many \emph{dissipative} weak solutions
$\vr$, $\vJ$ in $(0,T) \times \Omega$,
\[
\vr(0, \cdot) = \vr_0,\ \vJ(0,\cdot) = \vJ_0.
\]

\eT

%}

{\bf Sketch of the proof.} 

\medskip

We start with the ansatz $\vr(0, \cdot) = \vr_0$, $\vc{U}(0, \cdot) = \vc{U}_0 \equiv 0$, which yields $Z \equiv 0$ in (\ref{c1+}), in particular, 
$\vr (t, \cdot) = \vr_0$ is the unique solution of (\ref{c1}) for any $t > 0$. As shown in Section \ref{c}, the technique of convex integration produces (infinitely many) weak solutions $[\tilde \vr \equiv \vr, \tilde {\vc{J}}]$ emanating from $\vr_0$, $\vc{U}_0$, satisfying the energy inequality (\ref{ei}) in the \emph{open} interval $(0,T)$, where, in addition, the rate of total dissipation may controlled by a suitable choice of the function $h$ in (\ref{ei}). As $\vr$ is constant in time, the desired initial
distribution of the momentum $\vc{J}_0$ can be taken as $\vc{J}_0 = \tilde{\vc{J}} (\tau, \cdot)$ at a suitable point $\tau \in (0,T)$ where the energy is continuous.
Finally, the method of convex integration can be applied to produce infinitely many solutions for the initial data $\vr_0$, $\vc{J}_0$ as soon as we show that the set of subsolutions $X_{0,\ep}$ is non-empty. However, with Lemma \ref{LDS} at hand, a suitable subsolution can be constructed exactly as in \cite[Section 4.2]{ChiFeiKre}.

\def\cprime{$'$} \def\ocirc#1{\ifmmode\setbox0=\hbox{$#1$}\dimen0=\ht0
  \advance\dimen0 by1pt\rlap{\hbox to\wd0{\hss\raise\dimen0
  \hbox{\hskip.2em$\scriptscriptstyle\circ$}\hss}}#1\else {\accent"17 #1}\fi}

%\bibliography{citace}
%\bibliographystyle{plain}

\end{document}